\documentclass[a4paper,UKenglish,cleveref]{lipics-v2021}

\usepackage{algorithm}
\usepackage{algpseudocode}
\algrenewcommand\algorithmicrequire{\textbf{Input:}}
\algrenewcommand\algorithmicensure{\textbf{Output:}}

\usepackage{placeins}
\usepackage{booktabs}
\usepackage{accents}
\newcommand{\ubar}[1]{\underaccent{\bar}{#1}}

\hideLIPIcs  

\title{Exact and Efficient Methods for Identifying Carcinogenic Multi-Hit Gene Combinations} 

\titlerunning{Exact and Efficient Methods for Identifying Carcinogenic Multi-Hit Gene Combinations} 

\author{Rick S. H. Willemsen}{Singapore University of Technology and Design, Engineering Systems and Design, Singapore}{rick_willemsen@sutd.edu.sg}{}{}

\author{Teresa L. F. Ho}{Institute of Molecular and Cell Biology (IMCB), Agency for Science, Technology and Research (A*STAR), Singapore}{ teresa_ho@a-star.edu.sg}{}{}

\author{Tenindra Abeywickrama}{RIKEN Center for Computational Science, Japan}{tenindra.abeywickrama@riken.jp}{}{}

\authorrunning{R.\,S.\,H. Willemsen, T.\,L.\,F. Ho and T. Abeywickrama} 

\Copyright{Rick S. H. Willemsen, Teresa L. F. Ho and Tenindra Abeywickrama} 



\category{} 

\relatedversion{} 

\nolinenumbers 

\EventEditors{John Q. Open and Joan R. Access}
\EventNoEds{2}
\EventLongTitle{42nd Conference on Very Important Topics (CVIT 2016)}
\EventShortTitle{CVIT 2016}
\EventAcronym{CVIT}
\EventYear{2016}
\EventDate{December 24--27, 2016}
\EventLocation{Little Whinging, United Kingdom}
\EventLogo{}
\SeriesVolume{42}
\ArticleNo{23}

\begin{document}

\maketitle

\begin{abstract}
Cancer is driven by an estimated two to nine gene mutations, known as multi-hit combinations. Due to the large number of genes, identifying these multi-hit combinations presents a computationally challenging problem, which has previously required supercomputers to enumerate all or most gene combinations. Recent work formulated this classification task as the Multi-Hit Cancer Driver Set Cover Problem (MHCDSCP) and solved it via column generation on a single commodity CPU, where a pricing problem identifies promising gene combinations \cite{willemsen_et_al:LIPIcs.CP.2026.57}. The main bottleneck is the pricing problem that has to be solved repeatedly. To solve the pricing problem to optimality, we propose a dedicated depth first search method that enumerates all possible gene combinations. Although enumerating gene combinations with depth first search have been performed on supercomputers with around 150,000 CPU cores \cite{prabhu2026looking}, our depth first search algorithm solves the pricing problem to optimality within seconds on a single CPU by leveraging dual information to prune large parts of the search space. By eliminating this bottleneck, our approach yields the first exact method for the MHCDSCP, which is guaranteed to find a global optimal solution.
\end{abstract}

\section{Introduction}
The multi-hit theory of carcinogenesis states that multiple mutated genes (hits) are required to cause cancer
\cite{armitage1954age, nordling1953new, sutherl1984multihit}, where the number of hits is estimated to be between two to nine
\cite{little2003stochastic, zhang2005estimating, tomasetti2015only, ashley1969two, anandakrishnan2019estimating}. Initially, identifying these multi-hit combinations was formulated as a weighted set cover problem
\cite{dash2019differentiating}. This has been solved using full or partial enumeration of all gene combinations, followed by a greedy heuristic
\cite{dash2019differentiating, al2020identifying, dash2021scaling, dash2023distributing, oles2025bigpicc, prabhu2026looking}. To solve realistic sized instances, these heuristic approaches require supercomputers. Recently, the problem was formulated as the Multi-Hit Cancer Driver Set Cover Problem (MHCDSCP) \cite{willemsen_et_al:LIPIcs.CP.2026.57}, for which several mixed integer programming (MIP) formulations and column generation based heuristics were provided. They can solve the problem to near-optimality on a single commodity CPU within a minute. 

In this paper, we accelerate the proposed column generation based methods by proposing a dedicated algorithm to solve the pricing problem to optimality. This is subsequently incorporated in an exact branch-and-price framework. This leads to the first exact method for MHCDSCP, meaning that global optimal solutions are guaranteed given sufficient computational time.

The remainder of this paper is structured as follows. In \Cref{section2a:literature}, we compare our method to other column generation based classification algorithms in the literature. In \Cref{section2:problemdescription} we present the notation and MIP formulation of the MHCDSCP, while in \Cref{section3:colgen} we present our branch-and-price framework. Computational results on real-world data are shown in \Cref{section4:results} and \Cref{section5:conclusion} concludes. 

\section{Related column generation methods for classification}\label{section2a:literature}

The binary classification problem considered in this paper is closely related to interpretable rule generation for classification. In this setting, each column represents a classification rule, while the pricing problem identifies promising rules to be added to a restricted master problem.

Prior work generated boolean decision rules to make binary classification interpretable \cite{dash2018boolean}. Their objective is to minimize classification loss while generating an interpretable rule set. The pricing problem is solved via heuristics, subsampling, and an MIP solver, while the restricted master problem is solved using an interior point method or MIP solver. 

This framework is extended to fair and interpretable binary classification \cite{lawless2020fair}. The objective is to minimize false negatives while penalizing the number of selected rules. As in \cite{dash2018boolean}, the pricing problem combines heuristics, subsampling and an MIP solver. To obtain an integral rule set, they subsequently solve an MIP restricted to at most 1,000 generated rules. 

The multiclass version with fair and interpretable rules is considered for large datasets \cite{rober2025rule}. Their objective is to minimize classification loss while penalizing the number of selected rules. Their pricing problem is solved heuristically and their approach is designed to obtain a weighted rule set rather than an integral solution. 

The identification of multi-hit gene combinations can be formulated as the MHCDSCP, which is a binary classification problem, where a gene combination can be seen as an interpretable rule \cite{willemsen_et_al:LIPIcs.CP.2026.57}. The objective is to maximize true positives while minimizing false positives. Similar to \cite{dash2018boolean}, the pricing problem is solved using a combination of heuristics, subsampling and an MIP solver. An integral solution is obtained by solving an MIP model containing all generated columns.

The approaches discussed above generally rely on heuristics, subsampling, or MIP solvers to solve the pricing problem. In contrast, we develop a dedicated depth first search algorithm for the MHCDSCP that solves the pricing problem to optimality. The time complexity can be expressed as $\mathcal{O}(|G|^{\bar{k}})$, where $|G|$ is the number of genes and $\bar{k}$ the maximum hit size. In the MHCDSCP, the maximum hit size is at most nine, and hence the pricing problem can be solved in polynomial time. A relevant complexity result is provided by \cite{rober2025rule}, who proved that the pricing problem in their multiclass classification framework is NP-hard, while it becomes polynomially solvable when the rule size is fixed. Although our pricing problem is different, their result highlights the importance of a bounded rule size, which is bounded by nine in the MHCDSCP.

Despite this polynomial-time bound, enumerating all or most gene combinations is known to require supercomputers, see for instance \cite{dash2019differentiating, al2020identifying, dash2021scaling, dash2023distributing, oles2025bigpicc}. Specifically, \cite{prabhu2026looking} perform enumeration using a depth first search algorithm with pruning and require around 150,000 CPU cores. However, we propose a dedicated depth first search method for the pricing problem that uses dual information to prune large parts of the search space. The resulting pricing problem method is both exact and fast. This exact pricing procedure allows us to embed the column generation approach in a branch-and-price algorithm that is guaranteed to find a global optimal solution to the MHCDSCP.

\section{Problem description}\label{section2:problemdescription}

\subsection{Notation}
Mathematically, the problem can be formulated as follows. Let $G$ be the set of genes and $S$ be the set of samples, which is partitioned into tumor and normal samples, $T$ and $N$, such that $S=T \cup N$ and $T\cap N=\emptyset$. The input to the problem can be represented by a binary matrix $\mathbf{M}$ of size $|S|\times |G|$, where entry $\mathbf{M}_{sg}=1$ indicates that sample $s$ has a mutated gene~$g$.

A gene combination is defined as a subset of genes. Let $C$ be the set of all possible gene combinations with cardinality between $\ubar{k}$ and $\bar{k}$. A gene combination $c\in C$ is said to cover a sample $s\in S$ if all genes in $c$ are mutated in $s$. For ease of notation, we define $C_s\subseteq C$ as the set of all gene combinations that cover sample $s$ and we define $S_c\subseteq S$ as the set of all samples covered by combination $c$.

\subsection{Multi-Hit Cancer Driver Set Cover Problem}

The aim in the MHCDSCP is to find a set of gene combinations $C^*\subseteq C$ that induces an accurate binary classification of tumor and normal samples. A sample is classified as a tumor if it is covered by at least one gene combination from set $C^*$ and as a normal sample otherwise. In the MHCDSCP we balance two criteria, namely maximizing the number of correctly classified tumor samples while minimizing the number of times a normal sample is misclassified.

We denote by $T^*\subseteq T$ and $N^*\subseteq N$ the sets of samples that are classified as tumor by the selected gene combinations $C^*$. The number of true positives, false negatives, false positives, and true negatives are given by $|T^*|$, $|T\setminus T^*|$, $|N^*|$, and $|N\setminus N^*|$, respectively. To evaluate the performance of a classification, we consider the performance measures in \Cref{section2.3:performancemeasures}.

\subsection{Mixed integer programming formulation}
The MHCDSCP can be modeled as a mixed integer program. Let $x_t$ be a binary variable denoting whether a tumor sample $t$ is correctly classified, while $y_n$ is an integer variable counting the number of times a normal sample $n$ is misclassified. The selected gene combinations are modeled using binary variables $z_c$. Furthermore, we limit the number of selected gene combinations by $\beta$, since a limited number of gene combinations already provide good solutions \cite{willemsen_et_al:LIPIcs.CP.2026.57}.

Since the performance
measures we consider are nonlinear, the objective is linearized by maximizing the number of correctly classified tumor samples (true positives)
and minimizing the number of times a normal sample is misclassified (which approximates false
positives). This can be modeled as the following mixed integer programming formulation \cite{willemsen_et_al:LIPIcs.CP.2026.57} 
\begin{align}
    \max\ & \sum_{t\in T} x_t - \sum_{n\in N} y_n ,  \label{form1:cons1} \\
    \text{s.t.}\ & x_t - \sum_{c\in C_t} z_c\leq 0, && \forall t\in T ,\label{form1:cons2} \\
    & \sum_{c\in C_n} z_c-  y_n   \leq 0 , && \forall n\in N, \label{form1:cons3} \\
    & \sum_{c\in C} z_c  \leq  \beta  ,\label{form1:cons4} \\
    & x_t \in  \{0,1\} ,&& \forall t\in T ,\label{form1:cons5} \\
    & y_n \in  \{0, 1, \dots, \beta\} ,&& \forall n\in N ,\label{form1:cons6} \\
    & z_c \in \{0,1\} , && \forall c\in C .\label{form1:cons7}
\end{align}
The objective \eqref{form1:cons1} maximizes the number of correctly classified tumor samples and minimizes the misclassification of normal samples. Constraint \eqref{form1:cons2} ensures that a tumor sample $t$ is classified as a tumor if it is covered by at least one gene combination, while \eqref{form1:cons3} counts the number of times a normal sample $n$ is misclassified. The number of selected gene combinations is limited to $\beta$ via \eqref{form1:cons4}. The domain of the variables is enforced via \eqref{form1:cons5}-\eqref{form1:cons7}.

Set $C$ contains all gene combinations of size $\ubar{k}$ up to $\bar{k}$, resulting in $\mathcal{O}(|G|^{\bar{k}})$ variables. Column generation is a well-known technique to solve mixed integer programs with a large number of variables. In the next section, we present a branch-and-price method, which is a column generation method that solves the original problem to global optimality. Although this problem can be formulated using other techniques, such as constraint programming, for this specific version of the problem integer programming formulation seems more suitable \cite{willemsen_et_al:LIPIcs.CP.2026.57}.

\subsection{Trivial benchmark}
A natural question when proposing computationally intensive methods that require supercomputers or advanced mathematical programming techniques like branch-and-price is how it compares with a simple baseline. To our knowledge, previous work has not reported a comparison with a {trivial} benchmark. In our computational results, we show that a trivial benchmark already attains performance comparable to the state-of-the-art while running within a few minutes on a single commodity CPU. 

The trivial benchmark consists of two steps. In the first step, we greedily generate a set of gene combinations $\tilde{C}\subseteq C$ using Algorithm 1 from \cite{willemsen_et_al:LIPIcs.CP.2026.57}. In the second step, we greedily select gene combinations that improve the objective function while $|C^*|\leq \beta$. This second step can be implemented as a doubly nested for-loop, which we speed up by heuristically pruning parts of the search space.

To speed up the computation of the nested for-loop we heuristically prune the search space. We sort the obtained gene combinations $\tilde{C}$ in decreasing order of their potential individual contribution, $T_c$, i.e., the number of tumor samples covered by gene combination $c$. Let $f$ denote the objective function, while $f_t$ counts the number of covered tumor samples. Let $C^*$ be the set of gene combinations corresponding to the best objective $f^*$. Let $\Delta_t$ be the current improvement in the number of covered samples when adding a gene combination $c$ and let $\Delta_t^*$ denote an estimated improvement bound. 

When adding a new gene combination $c$ to $C^*$, the best possible improvement in the objective function is if all $T_c$ samples are covered due to $c$ only. If $T_c$ drops below the bound $\Delta_t^*$, we prune and stop the search.

\begin{algorithm}[h]
\caption{Trivial benchmark}
\label{alg:trivialBenchmark}
\begin{algorithmic}[1]
\State greedily generate gene combinations $\tilde{C}$ using Algorithm 1 in \cite{willemsen_et_al:LIPIcs.CP.2026.57}
\State initialize $C^*\leftarrow \emptyset$, $f^* \leftarrow -\infty$, $f_t^* \leftarrow -\infty$ and $\Delta_t^* \leftarrow -\infty$
\State sort $\tilde{C}$ in decreasing order of $T_c$
\While{$|C^*|<\beta$}
\State $c^*\leftarrow \emptyset$
    \For{$c\in \tilde{C}$}
    \If{$T_c < \Delta_t^*$} 
    \State break
    \EndIf
    \If{$f( C^*\cup c)>f^*$}
    \State $c^* \leftarrow  c$
        \State $f^*\leftarrow f(C^*\cup c^*)$
    \EndIf
    \State $\Delta_t\leftarrow  f_t(C^* \cup c)- f_t^*$
    \If{$\Delta_t> \Delta_t^*$}
    \State $\Delta_t^*\leftarrow \Delta_t$
    \EndIf
\EndFor
    \If{$c^*=\emptyset$}
    \State break
    \EndIf
\State $C^*\leftarrow C^* \cup c^*$
\State $f^* \leftarrow f(C^*) $
\State $f_t^* \leftarrow f_t(C^*)$
\EndWhile

\end{algorithmic}
\end{algorithm}
\FloatBarrier

\section{Branch-and-price}\label{section3:colgen}

\subsection{Overview}
Instead of solving the mixed integer program \eqref{form1:cons1}-\eqref{form1:cons7}, we solve a so-called restricted master problem on a subset of gene combinations $\tilde{C}\subseteq C$. Solving the restricted master problem provides information on each constraint in terms of dual values. These dual values can be used in a so-called pricing problem to generate gene combinations not yet in $\tilde{C}$ that may improve the objective value. Such a new gene combination is called a column and is added to $\tilde{C}$ in the restricted master problem. If the pricing problem is solved to optimality, iteratively repeating this procedure yields an optimal (possibly fractional) solution to the linear programming relaxation and provides a valid upper bound to the optimal objective value. By incorporating this procedure in a branch-and-bound framework, one obtains branch-and-price, which is guaranteed to provide a global optimal solution given enough time.

\subsection{Restricted master problem}
The restricted master problem is obtained by taking the linear relaxation of the mixed integer program restricted to a set of gene combinations $\tilde{C}$. In addition, we relax constraint \eqref{form1:cons7} to $z_c\geq 0$ noting that $z_c$ will not be set larger than one in an optimal solution. This leads to the following formulation for the restricted master problem \cite{willemsen_et_al:LIPIcs.CP.2026.57}
\begin{align}
    \max\ & \sum_{t\in T} x_t  - \sum_{n\in N} y_n , \label{form2:cons1} \\
    \text{s.t.}\ & x_t - \sum_{c\in \tilde{C}_t} z_c \leq 0, && \forall t\in T ,\label{form2:cons2} \\
    & \sum_{c\in \tilde{C}_n} z_c -  y_n \leq 0 , && \forall n\in N, \label{form2:cons3} \\
    & \sum_{c\in \tilde{C}} z_c  \leq  \beta  ,\label{form2:cons4} \\
    & x_t \in  [0, 1], && \forall t\in T ,\label{form2:cons5} \\
    & y_n \geq 0 ,&& \forall n\in N ,\label{form2:cons6} \\
    & z_c \geq 0 , && \forall c\in \tilde{C} .\label{form2:cons7} 
\end{align}

\subsection{Pricing problem formulation for the identification of gene combinations}
In the pricing problem, we aim to find a gene combination of size $\ubar{k}$ to $\bar{k}$ with so-called positive reduced cost that can be added to the restricted master problem or prove that no such gene combination exists. Solving the pricing problem to optimality guarantees that a gene combination with positive reduced cost is found if one exists.

Let $\pi_t\geq 0$, $\mu_n \geq 0$ and $\lambda \geq 0$ be the dual variables corresponding to \eqref{form2:cons2}-\eqref{form2:cons4}. The reduced cost of a gene combination $c$ can be calculated as
\begin{align*}
   RC(c) = \sum_{t\in T_c} \pi_t - \sum_{n\in N_c} \mu_n - \lambda.
\end{align*}
Note that $T_c$ and $N_c$ are the tumor and normal samples covered by gene combination $c$.

Let $u_g$ be a binary decision variable indicating that we select gene $g$ in the gene combination. In addition, let $v_t$ and $v_n$ be binary decision variables indicating that tumor sample $t$ or normal sample $n$ is covered by the constructed gene combination. The pricing problem can be formulated as a mixed integer program \cite{willemsen_et_al:LIPIcs.CP.2026.57} 
\begin{align}
    \max_{}\ & \sum_{t\in T} \pi_t v_t - \sum_{n\in N} \mu_n v_n - \lambda ,\label{formpp:cons1} \\
    \text{s.t.}\ &   \sum_{g\in G} u_g \geq \ubar{k}, \label{formpp:cons2} \\
    &  \sum_{g\in G} u_g \leq \bar{k} ,\label{formpp:cons3} \\
    & v_t \leq 1 - \frac{1}{\bar{k}} \sum_{g\in G: \mathbf{M}_{tg}=0} u_g , && \forall t\in T,\label{formpp:cons4} \\
    & v_n \geq 1 - \sum_{g\in G: \mathbf{M}_{ng}=0} u_g , && \forall n\in N,\label{formpp:cons6} \\
    & v_t \in \{0,1\} ,&&\forall t\in T,\label{formpp:cons7}  \\
    & v_n \in \{0,1\} ,&&\forall n\in N,\label{formpp:cons8}  \\
    & u_g \in \{0,1\} ,&&\forall g\in G.\label{formpp:cons9} 
\end{align}
Recall that $M_{sg}=1$ indicates that sample $s$ has a mutated gene $g$.

The objective \eqref{formpp:cons1} maximizes the reduced cost of the selected gene combination. Constraints \eqref{formpp:cons2}-\eqref{formpp:cons3} ensure that the selected gene combination contains between $\ubar{k}$ and $\bar{k}$ genes. A tumor sample is covered only if it contains all genes in the selected gene combination, which is modeled by \eqref{formpp:cons4}. On the other hand, \eqref{formpp:cons6} enforces that a normal sample is covered if it contains all genes from the constructed gene combination. The domain of the variables is modeled by \eqref{formpp:cons7}-\eqref{formpp:cons9}.

\subsection{Solving the pricing problem using depth first search}
The pricing problem can be solved first heuristically via a variable neighborhood search heuristic after which optimality should be proven by solving it with a mixed integer programming solver \cite{willemsen_et_al:LIPIcs.CP.2026.57}. In this manuscript, we propose to solve the pricing problem in polynomial time using an exact dedicated algorithm, namely depth first search (DFS). For a fixed ${k}$, enumeration using depth first search has complexity $\mathcal{O} (|G|^{{k}})$, which is polynomial. Note that enumeration of all gene combinations has been studied in the literature and requires supercomputers for large realistic instances \cite{dash2019differentiating, al2020identifying, dash2021scaling, dash2023distributing, oles2025bigpicc, prabhu2026looking}. However, we can use dual information to prune large parts of the search space, thereby making it feasible to solve the pricing problem using depth first search enumeration.

We obtain an upper bound on the reduced cost as follows
\begin{align*}
   RC(c) &= \sum_{t\in T_c} \pi_t v_t - \sum_{n\in N_c} \mu_n v_n - \lambda ,\\
  &\leq   \sum_{t\in T_c} \pi_t v_t - \lambda,\\
  &= \overline{RC}(c) .
\end{align*}
Adding genes to a gene combination can only remove covered tumor samples, since $\pi_t\geq 0$ the reduced cost upper bound $\overline{RC}(c)$ is monotonically decreasing as genes are added. We can use this in the DFS to prune that branch of the search space if $\overline{RC}(c)\leq LB^*$, where $LB^*$ is the reduced cost of the best solution found so far. The DFS algorithm is summarized in \Cref{alg:DFS}.

Intuitively, the main reasons why DFS can be fast are twofold. First, the maximum depth of the search tree equals the hit size $k$, which is at most nine in practice. Second, the pruning rule can easily prune large parts of the search space. Additional speed up is obtained by sorting the genes in descending order of $RC(\{g\})$.

\begin{algorithm}
\caption{Depth first search (DFS)}
\label{alg:DFS}
\begin{algorithmic}[1]
\Require start index $i$, gene combination $c$, hit size $k$ and an ordered list of genes ${G}$

\If{$\overline{RC}(c)\leq LB^*$}
\State return
\EndIf

\If{$|c|=k$}
\If{${RC}(c) > LB^*$}
\State a new best gene combination $c$
\State $LB^* = {RC}(c)$
\EndIf
\State return
\EndIf

\For{$j \leftarrow i$ to $|{G}|$}
    \State $g_j\in {G}$
    \State DFS($j+1$, $c \cup \{g_j\}$, $k$, ${G}$)
\EndFor

\end{algorithmic}
\end{algorithm}
\FloatBarrier

For better performance, we also perform subsampling. We first solve DFS with a subsampled list $\tilde{G}\subseteq G$ and for a fixed hit size, we loop over all possible hit sizes from $\ubar{k}$ up to $\bar{k}$. If no positive reduced cost is found, then we repeat it over the entire list of genes ${G}$. The pricing problem procedure is explained in \Cref{alg:PP}. Similar to \cite{willemsen_et_al:LIPIcs.CP.2026.57}, we initialize $G_{\text{VNS}}$ with initial genes and we solve the pricing problem using a variable neighborhood search heuristic on the restricted gene set $G_{\text{VNS}}\cup G_r$, where $G_r$ is initially empty. Every time we find a gene combination we add all these genes to $G_r$. If no positive reduced cost gene combination can be identified, we apply DFS on the restricted set. If still no gene combination with positive reduced cost is found we solve DFS on the entire gene set $G$ in order to solve the pricing problem to optimality. This makes the proposed pricing problem procedure exact.

\begin{algorithm}
\caption{Exact pricing problem procedure}
\label{alg:PP}
\begin{algorithmic}[1]
\Require Gene sets $G_{r}$ and $G_{\text{VNS}}$
\Ensure A gene combination $c$ with positive reduced cost

\State $c\gets$ apply variable neighborhood search restricted to $G_{\text{VNS}}\cup G_{r}$ (Algorithm 5 in \cite{willemsen_et_al:LIPIcs.CP.2026.57})
\Statex \hspace{\algorithmicindent} \textbf{if} $c \neq \emptyset$ \textbf{goto} step 4.
\State $c\gets$ apply DFS restricted to $G_{r}$
\Statex \hspace{\algorithmicindent} \textbf{if} $c \neq \emptyset$ \textbf{goto} step 4.
\State $c\gets$ apply DFS using entire gene set $G$
\State \Return $c$
\end{algorithmic}
\end{algorithm}
\FloatBarrier

Note that the MIP solver for the pricing problem in \cite{willemsen_et_al:LIPIcs.CP.2026.57} returns the first positive reduced cost column. In contrast, we return the most positive reduced cost column in line 3 of \Cref{alg:PP}. In addition, we can return multiple positive reduced cost columns by saving intermediate solutions.

\subsection{Branching rule}
When all gene combination variables $z_c$ in the mixed integer program are fixed, the continuous $x_t$ and $y_n$ variables take integral values in an optimal solution \cite{willemsen_et_al:LIPIcs.CP.2026.57}. Therefore, it suffices to branch on the $z_c$ variables.

In this paper, we identify the most fractional variable corresponding to a gene combination $c$ and create two branches, namely $z_c=0$ and $z_c=1$. We first branch on the enforced variable $z_c=1$. We have to ensure that the pricing problem is correct when adding these constraints to the restricted master problem. When $z_c=1$, no changes have to be made, since the corresponding gene combination $c$ cannot have positive reduced cost. When $z_c=0$, we have to make an adjustment to the pricing problem. Specifically, we have to ensure that gene combination $c$ cannot be generated. 

Note that if one were to solve the pricing problem using an MIP solver, no-good cuts are needed, which would slow down the method considerably. In DFS, before evaluating a candidate combination, we check whether it is forbidden by any branching constraint $z_c=0$ and skip it if so. We perform this check before the recursive call on line 13 of \Cref{alg:DFS}.

\section{Computational results}\label{section4:results}

\subsection{Experimental setup}
We evaluate our proposed methods on real-world data obtained from \cite{oles2025bigpicc}, which is publicly available via the Oak Ridge National Laboratory Gitlab\footnote{\url{https://code.ornl.gov/vo0/bigpicc}}. The dataset contains two cancer types, whose characteristics are summarized in \Cref{tab:summarystatistics}.

\begin{table}[h]
\caption{Overview of characteristics of the dataset. We report the hit sizes ($\ubar{k}$–$\bar{k}$), number of genes, samples and mutations for each cancer type.}
\label{tab:summarystatistics}
\centering
\begin{tabular}{llrrr}
\toprule
 cancer type & \#hits        & \#genes & \#samples & \#mutations \\ \midrule
 BRCA        & 3, 2–3, 6–7 & 19,411 & 1242    & 837,207    \\
 PRAD        & 6–7, 7      & 17,274 & 752     & 476,115    \\ 
\bottomrule
\end{tabular}
\end{table}
\FloatBarrier

We compare our methods against two heuristics that can be run on a single commodity CPU \cite{willemsen_et_al:LIPIcs.CP.2026.57}, for which the code can be found online\footnote{\url{https://github.com/rickwillemsen/MHCDSCP}}. Their fastest method is an MIP heuristic while the best performing method based on performance measures is a root node heuristic. Both the MIP heuristic and trivial benchmark require an initial set of gene combinations. We initialize both methods with the same $100,000$ greedily generated gene combinations. For completeness we also report the results from a recent method named BiGPICC, which was executed on 280 nodes on the Summit supercomputer \cite{oles2025bigpicc}.

In all our experiments, we limit the number of selected gene combinations to $\beta=10$ following \cite{willemsen_et_al:LIPIcs.CP.2026.57}. Unless specified otherwise, the training set contains 75\% of the data while the remaining 25\% is used as test set. All experiments were conducted on a single 12-core Intel Core Ultra 5 225U processor running at 1.50 GHz and 32 GB RAM. All
MIP formulations are solved using the commercial solver CPLEX 20.10.

\subsection{Trivial benchmark}
For fairness, we run the MIP heuristic, root node heuristic and trivial benchmark on the same training-test split. Following \cite{willemsen_et_al:LIPIcs.CP.2026.57}, the MIP heuristic has 30 second time limit for solving the model, while the root node heuristic has a 300 second time limit. The trivial benchmark is run without a time limit. For BiGPICC we show the results as reported in \cite{oles2025bigpicc}, although they also use 75\%-25\% training-test split the exact split is different. 

The results on the test set performance are shown in \Cref{tab:trivialBenchmark}, where we report the test performance measured by the Matthews correlation coefficient (MCC)~\cite{chicco2020advantages}. The trivial benchmark performs better in three instances compared to BiGPICC, while the average test set performance is slightly lower. Compared to the MIP heuristic the trivial benchmark performs twice strictly better and once has the same performance, also here the average MCC value is quite comparable. The trivial benchmark also performs twice better than the root node heuristic, but the root node heuristic has a better average test set performance. 

When looking at the computation time, the trivial benchmark is the fastest on average. While the more advanced root node heuristic obtains the best solutions, it is considerably slower.

\begin{table}[h]
\centering
\caption{Comparison of test set performance between BiGPICC (results reported from \cite{oles2025bigpicc}), MIP heuristic \cite{willemsen_et_al:LIPIcs.CP.2026.57}, root node heuristic \cite{willemsen_et_al:LIPIcs.CP.2026.57}, and our proposed trivial benchmark}
\label{tab:trivialBenchmark}
\resizebox{\textwidth}{!}{%
\begin{tabular}{llrrrrrrrrrrrrrrr}
\toprule
\shortstack[l]{cancer\\ type} & \shortstack[l]{hit\\ range} &  & \multicolumn{2}{c}{\shortstack[c]{BiGPICC\\ (results reported from \cite{oles2025bigpicc})}} &  & \multicolumn{3}{c}{MIP heuristic \cite{willemsen_et_al:LIPIcs.CP.2026.57}} &  & \multicolumn{3}{c}{root node heuristic \cite{willemsen_et_al:LIPIcs.CP.2026.57}}  &  & \multicolumn{3}{c}{trivial benchmark} \\ \cline{4-5} \cline{7-9} \cline{11-13} \cline{15-17} 
            &           &  & \#comb        & MCC         &  & \#comb  & MCC    & time &  & \#comb   & MCC      & time     &  & \#comb       & MCC         & time       \\ \midrule
\multirow{3}{*}{BRCA}            & 2–3       &  & 5             & 0.900       &  & 10      & 0.797  & 37.9 &  &     10     &  0.896        &     30.6     &  & 8            & 0.917       & 73.0       \\
        & 3         &  & 9             & 0.900       &  & 9       & 0.917  & 39.7 &  &    10      &  0.909        &   39.6       &  & 6            & 0.917       & 41.5       \\
            & 6–7       &  & 34            & 0.749       &  & 10      & 0.560  & 34.0 &  &      10    & 0.720         &   300.0       &  & 10           & 0.603       & 14.7       \\ \midrule
\multirow{2}{*}{PRAD}        & 6–7       &  & 45            & 0.436       &  & 10      & 0.642  & 23.3 &  &       10   &    0.571      &  300.0        &  & 10           & 0.511       & 7.7        \\
            & 7         &  & 62            & 0.487       &  & 10      & 0.564  & 30.6 &  &     10     &     0.692     &  300.0        &  & 10           & 0.390       & 4.8        \\ \midrule
avg         &           &  & 31.0          & 0.694       &  & 9.8     & 0.696  & 33.1 &  & 10.0 & 0.758 & 194.1 &  & 8.8          & 0.667       & 28.3      \\ \bottomrule
\end{tabular}
}
\end{table}
\FloatBarrier

\subsection{Analysing the pricing problem}

The root node heuristic in \cite{willemsen_et_al:LIPIcs.CP.2026.57} uses an MIP solver to solve the pricing problem, where it aims to find a gene combination with positive reduced cost. In contrast, our DFS method finds an optimal reduced cost gene combination. Therefore, to evaluate the impact of our proposed pricing procedure, we use the same root node heuristic with a time limit of 300 seconds and only replace the pricing problem solver. We compare the following four strategies: MIP(1), vanilla DFS (1), DFS(1), DFS(20), where the number in brackets denotes the number of positive reduced cost gene combinations that are returned. Vanilla DFS means that it performs DFS without sorting the genes. We compare the methods in terms of both computation time and the quality of the bounds, measured by best feasible solution (obj) and current LP value. 

The computational results are shown in \Cref{tab:pp}. We can observe that a vanilla DFS performs worse than \cite{willemsen_et_al:LIPIcs.CP.2026.57}, however, when sorting the genes, DFS(1), we are quickly able to find a lower bound that helps prune the search space. In addition, by adding 20 instead of 1 column per iteration we obtain more improvements in the objective. For DFS(20), the average solution time in seconds per pricing problem iteration ranges from 0.04 to 5.2 seconds. These results show that the pricing problem can now be solved efficiently within seconds on real-world instances. Having a fast pricing procedure makes a branch-and-price method viable.

\begin{table}[h]
\centering
\caption{Root node heuristic results when changing the pricing problem solver. We report the best feasible objective value, current LP value (bound) and computation time for four pricing problem configurations.}
\label{tab:pp}
\resizebox{\textwidth}{!}{%
\begin{tabular}{llrrrrrrrrrrrrrrrr}
\toprule
cancer type & hit range &  & \multicolumn{3}{c}{MIP(1) \cite{willemsen_et_al:LIPIcs.CP.2026.57}} &  & \multicolumn{3}{c}{vanilla DFS(1)} &  & \multicolumn{3}{c}{DFS(1)} &  & \multicolumn{3}{c}{DFS(20)} \\ \cline{4-6} \cline{8-10} \cline{12-14} \cline{16-18} 
            &           &  & obj    & bound   & time    &  & obj      & bound      & time       &  & obj    & bound   & time    &  & obj    & bound    & time    \\ \midrule
\multirow{3}{*}{BRCA}        & 2–3        &  & 677    & 683.0   & 30.6    &  & 677      & 683.0      & 15.8       &  & 677    & 683.0   & 7.3     &  & 677    & 683.0    & 7.1     \\
            & 3         &  & 675    & 683.0   & 39.6    &  & 676      & 683.0      & 8.1        &  & 676    & 683.0   & 9.2     &  & 676    & 683.0    & 8.5     \\
            & 6–7        &  & 649    & 664.5   & 300.0   &  & 645      & 662.8      & 300.0      &  & 655    & 675.1   & 300.0   &  & 657    & 677.1    & 300.0   \\ \midrule
\multirow{2}{*}{PRAD}        & 6–7        &  & 284    & 293.7   & 300.0   &  & 272      & 274.1      & 300.0      &  & 281    & 293.4   & 300.0   &  & 284    & 296.9    & 300.0   \\
            & 7         &  & 265    & 270.4   & 300.0   &  & 262      & 269.4      & 300.0      &  & 271    & 283.3   & 300.0   &  & 280    & 290.9    & 300.0   \\ \midrule
avg         &           &  &        &         & 194.1   &  &          &            & 184.8      &  &        &         & 183.3   &  &        &          & 183.1  \\ \bottomrule
\end{tabular}
}
\end{table}
\FloatBarrier

\subsection{Branch-and-price}
Finally, we apply branch-and-price on a smaller training sample, ranging from 30\% to 50\% of the data for a few selected cancer types and hit ranges. The computational results can be found in \Cref{tab:branchAndPrice} where the branch-and-price procedure has a 300 second time limit. These results show that some moderate sized instances can be solved to optimality by branch-and-price within the time limit. For the instances solved to optimality, the number of branches is strictly positive, indicating that solving the root node alone is insufficient to certify optimality. The largest instance solved to optimality contains 496 samples. 

These results indicate that branch-and-price can solve moderate sized instances to optimality and obtain valid bounds on larger instances within the time limit. This makes the approach potentially useful as a subroutine in applications such as cross validation, where multiple subsets of the data need to be solved.

\begin{table}[h]
\centering
\caption{Computational results when applying branch-and-price on a subset of the training data.}
\label{tab:branchAndPrice}
\begin{tabular}{llrrrrrrrr}
\toprule
cancer type     & hit range     &  & train (\%) & sample &  & obj & bound & time & \#branches \\ \midrule
\multirow{4}{*}{BRCA} & 2–3 &  & 40         & 496    &  & 455       & 455         & 7.3             & 27                \\
     & 2–3    &  & 50         & 620    &  & 544       & 547         & 300.0           & 271               \\ \cline{2-10}
     & 6–7 &  & 30         & 372    &  & 273       & 273         & 260.0           & 397               \\
     &  6–7   &  & 40         & 496    &  & 357       & 364         & 300.0           & 21                \\ \midrule
\multirow{2}{*}{PRAD} & 6–7 &  & 30         & 225    &  & 126       & 126         & 76.4            & 83                \\
     & 6–7    &  & 40         & 300    &  & 166       & 168         & 300.0           & 32       \\ \bottomrule         
\end{tabular}
\end{table}
\FloatBarrier

\section{Conclusion}\label{section5:conclusion}
In this paper, we introduce two new solution methods for the MHCDSCP, a trivial benchmark and a branch-and-price method.

The trivial benchmark attains test set performance comparable to heuristic solvers that require supercomputers or specialized software. This makes the trivial benchmark straightforward for practitioners to run locally and provides an accessible baseline for further research. In addition, the benchmark can serve as a starting point for studying more complex variants of the MHCDSCP that cannot be handled by a similar trivial procedure. For instance, some nonlinear constraints can be addressed using constraint programming \cite{willemsen_et_al:LIPIcs.CP.2026.57}. 

The branch-and-price method can successfully solve moderate sized instances to optimality within a few minutes. This is enabled due to a dedicated depth first search algorithm that can solve the pricing problem to optimality by exploiting dual information. An exact method provides a tool for testing hypotheses, such as the effect of the assumed hit size.



\bibliography{lipics-v2021-sample-article}

\appendix

\section{Performance measures}\label{section2.3:performancemeasures}
We consider the performance measures specificity, sensitivity, F1 score, and the Matthews correlation coefficient (MCC)~\cite{chicco2020advantages}, defined as
\begin{align*}
    \text{sensitivity} &= \frac{\text{TP}}{\text{TP} + \text{FN}},\\
    \text{specificity} &= \frac{\text{TN}}{\text{TN} + \text{FP}},\\
    \text{precision} &= \frac{\text{TP}}{\text{TP} + \text{FP}},\\
    \text{F1 score} &= 2\cdot \frac{\text{precision} \cdot \text{sensitivity}}{\text{precision} + \text{sensitivity}},\\
    \text{MCC} &= \frac{\text{TP} \cdot \text{TN} - \text{FP} \cdot \text{FN}}{\sqrt{ (\text{TP} + \text{FP})(\text{TP} + \text{FN})(\text{TN} + \text{FP})(\text{TN} + \text{FN})}},
\end{align*}
where true positives, false negatives, false positives, and true negatives are denoted by TP, FN, FP, and TN, respectively.

\end{document}